# HITTING TIMES FOR GAUSSIAN PROCESSES

By Laurent Decreusefond and David Nualart[1]

*GET/Telecom Paris and University of Kansas*

We establish a general formula for the Laplace transform of the hitting times of a Gaussian process. Some consequences are derived, and particular cases like the fractional Brownian motion are discussed.

**1. Introduction.** Consider a zero mean continuous Gaussian process $(X_t, t \geq 0)$, and for any $a > 0$, we denote by $\tau_a$ the hitting time of the level $a$ defined by

$$(1.1) \qquad \tau_a = \inf\{t \geq 0 : X_t = a\} = \inf\{t \geq 0 : X_t \geq a\}.$$

Thus, the map $(a \mapsto \tau_a)$ is left-continuous and increasing, hence, with right limits. The map $(a \mapsto \tau_{a^+})$ is right continuous where

$$\tau_{a^+} = \lim_{b \downarrow a, b > a} \tau_a = \inf\{t \geq 0 : X_t > a\}.$$

Little is known about the distribution of $\tau_a$. It is explicitly known in particular cases like the Brownian motion. If $X$ is a fractional Brownian motion with Hurst parameter $H$, there is a result by Molchan [5] which stands that

$$P(\tau_a > t) = t^{-(1-H)+o(1)}$$

as $t$ goes to infinity.

When $X$ is a standard Brownian motion, it is well known that

$$(1.2) \qquad E(\exp(-\alpha \tau_a)) = \exp(-a\sqrt{2\alpha})$$

for all $\alpha > 0$. This result is easily proved using the exponential martingale

$$M_t = \exp(\lambda B_t - \tfrac{1}{2}\lambda^2 t).$$

Received September 2006; revised January 2007.
[1]Supported in part by the NSF Grant DMS-06-04207.

By Doob's optional stopping theorem applied at time $t \wedge \tau_a$ and letting $t \to \infty$, one gets $1 = E(M_{\tau_a}) = E(\exp(\lambda B_{\tau_a} - \lambda^2 \tau_a/2))$. Since $B_{\tau_a} = a$, we thus obtain (1.2). If we consider a general Gaussian process $X_t$, the exponential process

$$M_t = \exp(\lambda X_t - \tfrac{1}{2}\lambda^2 V_t),$$

where $V_t = E(X_t^2)$ is no longer a martingale. However, it is equal to 1 plus a divergence integral in the sense of Malliavin calculus. The aim of this paper is to take advantage of this fact in order to derive a formula for $E(\exp(-\tfrac{1}{2}\lambda^2 V_{\tau_a}))$. We derive an equation involving this expectation in Theorem 3.4, under rather general conditions on the covariance of the process. As a consequence, we show that if the partial derivative of the covariance is nonnegative, then $E(\exp(-\tfrac{1}{2}\lambda^2 V_{\tau_a})) \leq 1$, which implies that $V_{\tau_a}$ has infinite moments of order $p$ for all $p \geq \tfrac{1}{2}$ and finite negative moments of all orders. In particular, for the fractional Brownian motion with Hurst parameter $H > \tfrac{1}{2}$, we have the inequality

$$E(\exp - \alpha \tau_a^{2H}) \leq \exp(-a\sqrt{2\alpha})$$

for all $\alpha, a > 0$.

The paper is organized as follows. In Section 2 we present some preliminaries on Malliavin calculus, and the main results are proved in Section 3.

**2. Preliminaries on Malliavin calculus.** Let $(X_t, t \geq 0)$ be a zero mean Gaussian process such that $X_0 = 0$ and with covariance function

$$R(s,t) = E(X_t X_s).$$

We denote by $\mathcal{E}$ the set of step functions on $[0, +\infty)$. Let $\mathcal{H}$ be the Hilbert space defined as the closure of $\mathcal{E}$ with respect to the scalar product

$$\langle \mathbf{1}_{[0,t]}, \mathbf{1}_{[0,s]} \rangle_{\mathcal{H}} = R(t,s).$$

The mapping $\mathbf{1}_{[0,t]} \longrightarrow X_t$ can be extended to an isometry between $\mathcal{H}$ and the Gaussian space $H_1(X)$ associated with $X$. We will denote this isometry by $\varphi \longrightarrow X(\varphi)$.

Let $\mathcal{S}$ be the set of smooth and cylindrical random variables of the form

(2.1) $$F = f(X(\phi_1), \ldots, X(\phi_n)),$$

where $n \geq 1$, $f \in \mathcal{C}_b^\infty(\mathbb{R}^n)$ ($f$ and all its partial derivatives are bounded), and $\phi_i \in \mathcal{H}$.

The *derivative operator* $D$ of a smooth and cylindrical random variable $F$ of the form (2.1) is defined as the $\mathcal{H}$-valued random variable

$$DF = \sum_{i=1}^n \frac{\partial f}{\partial x_i}(X(\phi_1), \ldots, X(\phi_n)) \phi_i.$$

The derivative operator $D$ is then a closable operator from $L^2(\Omega)$ into $L^2(\Omega; \mathcal{H})$. The Sobolev space $\mathbb{D}^{1,2}$ is the closure of $\mathcal{S}$ with respect to the norm

$$\|F\|_{1,2}^2 = E(F^2) + E(\|DF\|_{\mathcal{H}}^2).$$

The *divergence operator* $\delta$ is the adjoint of the derivative operator. We say that a random variable $u$ in $L^2(\Omega; \mathcal{H})$ belongs to the domain of the divergence operator, denoted by Dom $\delta$, if

$$|E(\langle DF, u \rangle_{\mathcal{H}})| \leq c_u \|F\|_{L^2(\Omega)}$$

for any $F \in \mathcal{S}$. In this case $\delta(u)$ is defined by the duality relationship

(2.2) $$E(F\delta(u)) = E(\langle DF, u \rangle_{\mathcal{H}}),$$

for any $F \in \mathbb{D}^{1,2}$.

Set $V_t = R(t,t)$. For any $\lambda > 0$, we define

$$M_t = \exp(\lambda X_t - \tfrac{1}{2}\lambda^2 V_t).$$

Formally, the Itô formula for the divergence integral, proved, for instance, in [1], implies that

(2.3) $$M_t = 1 + \lambda \delta(M \mathbf{1}_{[0,t]}),$$

where $M\mathbf{1}_{[0,t]}$ represents the process $(s \mapsto M_s \mathbf{1}_{[0,t]}(s), s \geq 0)$. However, the process $M\mathbf{1}_{[0,t]}$ does not belong, in general, to the domain of the divergence operator. This happens, for instance, in the following basic example.

EXAMPLE 1. Fractional Brownian motion with Hurst parameter $H \in (0,1)$ is a zero mean Gaussian process $(B_t^H, t \geq 0)$ with the covariance

(2.4) $$R_H(t,s) = \tfrac{1}{2}(t^{2H} + s^{2H} - |t-s|^{2H}).$$

In this case, the processes $(B_s^H \mathbf{1}_{[0,t]}(s), s \geq 0)$ and $(\exp(\lambda B_s^H - \tfrac{1}{2}\lambda^2 s^{2H})\mathbf{1}_{[0,t]}(s), s \geq 0)$ do not belong to $L^2(\Omega; \mathcal{H})$ if $H \leq \tfrac{1}{4}$ (see [2]).

In order to define the divergence of $M\mathbf{1}_{[0,t]}$ and to establish formula (2.3), we introduce the following additional property on the covariance function of the process $X$.

(H0) The covariance function $R(t,s)$ is continuous, the partial derivative $\frac{\partial R}{\partial s}(s,t)$ exists in the region $\{0 < s, t, s \neq t\}$, and for all $T > 0$,

$$\sup_{t \in [0,T]} \int_0^T \left|\frac{\partial R}{\partial s}(s,t)\right| ds < \infty.$$

Notice that this property is satisfied by the covariance (2.4) for all $H \in (0,1)$.

Define

(2.5) $$\delta_t M = \frac{1}{\lambda}(M_t - 1).$$

The following proposition asserts that $\delta_t M$ satisfies an integration by parts formula, and in this sense, it coincides with an extension of the divergence of $M\mathbf{1}_{[0,t]}$.

PROPOSITION 2.1. *Suppose that* (H0) *holds. Then, for any $t > 0$, and for any smooth and cylindrical random variable of the form $F = f(X_{t_1}, \ldots, X_{t_n})$, we have*

(2.6) $$E(F\delta_t M) = E\left(\sum_{i=1}^n \frac{\partial f}{\partial x_i}(X_{t_1}, \ldots, X_{t_n}) \int_0^t M_s \frac{\partial R}{\partial s}(s, t_i)\, ds\right).$$

PROOF. First notice that condition (H0) implies that the right-hand side of equation (2.6) is well defined. Then, it suffices to show equation (2.6) for a function of the form

$$f(x_1, \ldots, x_n) = \exp\left(\sum_{i=1}^n \lambda_i x_i\right),$$

where $\lambda_i \in \mathbb{R}$. In this case we have, for all $0 < t_1 < \cdots < t_n$,

$$\frac{1}{\lambda} E(F(M_t - 1))$$

$$= \frac{1}{\lambda} \exp\left\{\frac{1}{2}\sum_{i=1}^n \lambda_i \lambda_j R(t_i, t_j)\right\}\left(\exp\left\{\sum_{i=1}^n \lambda \lambda_i R(t, t_i)\right\} - 1\right)$$

$$= \sum_{i=1}^n \int_0^t \exp\left\{\frac{1}{2}\sum_{i=1}^n \lambda_i \lambda_j R(t_i, t_j) + \lambda \sum_{i=1}^n \lambda_i R(s, t_i)\right\} \lambda_i \frac{\partial R}{\partial s}(s, t_i)\, ds$$

$$= \int_0^t E\left(\sum_{i=1}^n \frac{\partial f}{\partial x_i}(X_{t_1}, \ldots, X_{t_n}) M_s \frac{\partial R}{\partial s}(s, t_i)\right) ds,$$

which completes the proof of the proposition. □

In many cases like in Example 1 with $H > \frac{1}{4}$, the process $M\mathbf{1}_{[0,t]}$ belongs to the space $L^2(\Omega; \mathcal{H})$, and then the right-hand side of equation (2.6) equals

$$E\langle DF, M\mathbf{1}_{[0,t]}\rangle_{\mathcal{H}}.$$

In this situation, taking into account the duality formula (2.2), equation (2.6) says that $M\mathbf{1}_{[0,t]}$ belongs to the domain of the divergence and $\delta(M\mathbf{1}_{[0,t]}) = \delta_t M$.

**3. Hitting times.** In this section we will assume the following conditions:

(H1) The partial derivative $\frac{\partial R}{\partial s}(s,t)$ exists and it is continuous in $[0,+\infty)^2$.
(H2) $\limsup_{t\to\infty} X_t = +\infty$ almost surely.
(H3) For any $0 \le s < t$, we have $E(|X_t - X_s|^2) > 0$.

Under these conditions, the process $X$ has a continuous version because

$$E(|X_t - X_s|^2) = R(t,t) + R(s,s) - 2R(s,t)$$
$$= \int_s^t \left[\frac{\partial R}{\partial u}(u,t) - \frac{\partial R}{\partial u}(u,s)\right] du$$
$$\le 2|t-s| \sup_{s \le u \le t}\left|\frac{\partial R}{\partial u}(u,t)\right|.$$

For any $a > 0$, we define the hitting time $\tau_a$ by (1.1). We know that $P(\tau_a < \infty) = 1$ by condition (H2). Set

(3.1) $$S_t = \sup_{s \in [0,t]} X_s.$$

From the results of [6], it follows that, for all $t > 0$, the random variable $S_t$ belongs to the space $\mathbb{D}^{1,2}$. Furthermore, condition (H3) allows us to compute the derivative of this random variable.

LEMMA 3.1. *For all $t > 0$, with probability one, the maximum of the process $X$ in the interval $[0,t]$ is attained in a unique point, that is, $\tau_{S_t} = \tau_{S_t^+}$ and $DS_t = \mathbf{1}_{[0,\tau_{S_t}]}$.*

PROOF. The fact that the maximum is attained in a unique point follows from condition (H3) and Lemma 2.6 in Kim and Pollard [4]. The formula for the derivative of $S_t$ follows easily by an approximation argument. □

We need the following regularization of the stopping time $\tau_a$. Suppose that $\varphi$ is a nonnegative smooth function with compact support in $(0, +\infty)$ and define for any $T > 0$

(3.2) $$Y = \int_0^\infty \varphi(a)(\tau_a \wedge T)\, da.$$

The next result states the differentiability of the random variable $Y$ in the sense of Malliavin calculus and provides an explicit formula for its derivative.

LEMMA 3.2. *The random variable $Y$ defined in* (3.2) *belongs to the space $\mathbb{D}^{1,2}$, and*

(3.3) $$D_r Y = -\int_0^{S_T} \varphi(y)\mathbf{1}_{[0,\tau_y]}(r)\, d\tau_y.$$

PROOF. Clearly, $Y$ is bounded. On the other hand, for any $r > 0$, we have

$$\{\tau_a > r\} = \{S_r < a\}.$$

Therefore, we can write using Fubini's theorem

$$Y = \int_0^\infty \varphi(a) \left( \int_0^{\tau_a \wedge T} d\theta \right) da = \int_0^T \left( \int_{S_\theta}^\infty \varphi(a) \, da \right) d\theta,$$

which implies that $Y \in \mathbb{D}^{1,2}$ because $S_\theta \in \mathbb{D}^{1,2}$, and

$$D_r Y = -\int_0^T \varphi(S_\theta) D_r S_\theta \, d\theta = -\int_0^T \varphi(S_\theta) \mathbf{1}_{[0, \tau_{S_\theta}]}(r) \, d\theta.$$

Finally, making the change of variable $S_\theta = y$ yields

$$D_r Y = -\int_0^{S_T} \varphi(y) \mathbf{1}_{[0, \tau_y]}(r) \, d\tau_y. \qquad \square$$

Notice that $M_Y = \exp(\lambda X_Y - \frac{1}{2}\lambda^2 V_Y)$. Hence, letting $t = Y$ in equation (2.5) and taking the mathematical expectation of both members of the equality yields

(3.4) $$E(M_Y) = 1 + \lambda E(\delta_t M|_{t=Y}).$$

We are going to show the following result which provides a formula for the left-hand side of equation (3.4).

LEMMA 3.3. *Assume conditions* (H1), (H2) *and* (H3). *Then, we have*

(3.5) $$E(M_Y) = 1 - \lambda E\left( M_Y \int_0^{S_T} \varphi(y) \frac{\partial R}{\partial s}(Y, \tau_y) \, d\tau_y \right).$$

PROOF. The proof will be done in two steps.
*Step* 1. We claim that for any function $p(x)$ in $\mathcal{C}_0^\infty(\mathbb{R})$ we have

(3.6) $$E(\delta_t M p(Y)) = -E\left( \int_0^t M_s p'(Y) \int_0^{S_T} \varphi(y) \frac{\partial R}{\partial s}(s, \tau_y) \, d\tau_y \, ds \right).$$

We can write $Y = \int_0^T \psi(S_\theta) \, d\theta$, where $\psi(x) = \int_x^\infty \varphi(a) \, da$. Consider an increasing sequence $D_n$ of finite subsets of $[0, T]$ such that their union is dense in $[0, T]$. Set $Y_n = \int_0^T \psi(S_\theta^n) \, d\theta$, and $S_\theta^n = \max\{X_t, t \in D_n \cap [0, \theta]\}$. Then, $Y_n$ is a Lipschitz function of $\{X_t, t \in D_n\}$. Hence, formula (2.6), which holds for Lipschitz functions, implies that

$$E(\delta_t M p(Y_n)) = -E\left( p'(Y_n) \int_0^T \varphi(S_\theta^n) \left( \int_0^t M_s \frac{\partial R}{\partial s}(s, \tau_{S_\theta^n}) \, ds \right) d\theta \right).$$

The function $r \to \int_0^t M_s \frac{\partial R}{\partial s}(s, r) \, ds$ is continuous and bounded by condition (H1). As a consequence, we can take the limit of the above expression as $n$ tends to infinity and we get

$$E(\delta_t M p(Y)) = -E\left(p'(Y) \int_0^T \varphi(S_\theta) \left(\int_0^t M_s \frac{\partial R}{\partial s}(s, \tau_{S_\theta}) \, ds\right) d\theta\right).$$

Finally, making the change of variable $S_\theta = y$ yields (3.6).

*Step* 2. We write

$$E(\delta_t M|_{t=Y}) = E\left(\lim_{\varepsilon \to 0} \int_{-\infty}^\infty \delta_t M p_\varepsilon(Y - t) \, dt\right),$$

where $p_\varepsilon(x)$ is an approximation of the identity, and by convention, we assume that $\delta_t M = 0$ if $t$ is negative. We can commute the expectation with the above limit by the dominated convergence theorem because

$$\int_{-\infty}^\infty |\delta_t M| p_\varepsilon(Y - t) \, dt = \int_{-\infty}^\infty \frac{1}{\lambda} |M_t - 1| p_\varepsilon(Y - t) \, dt$$

$$\leq \frac{1}{\lambda} \sup_{0 \leq t \leq T+1} (|M_t| + 1),$$

if the support of $p_\varepsilon(x)$ is included in $[-\varepsilon, \varepsilon]$, and $\varepsilon \leq 1$. Hence,

(3.7) $$E(\delta_t M|_{t=Y}) = \lim_{\varepsilon \to 0} \int_{-\infty}^\infty E(\delta_t M p_\varepsilon(Y - t)) \, dt.$$

Using formula (3.6) yields

(3.8)
$$E(\delta_t M p_\varepsilon(Y - t))$$
$$= -\int_0^t E\left(p'_\varepsilon(Y - t) M_s \left(\int_0^{S_T} \varphi(y) \frac{\partial R}{\partial s}(s, \tau_y) \, d\tau_y\right)\right) ds.$$

Hence, substituting (3.8) into (3.7) and integrating by parts, we obtain

$$E(\delta_t M|_{t=Y})$$
$$= -\lim_{\varepsilon \to 0} E\left(\int_{-\infty}^\infty p'_\varepsilon(Y - t) \left(\int_0^t M_s \left(\int_0^{S_T} \varphi(y) \frac{\partial R}{\partial s}(s, \tau_y) \, d\tau_y\right) ds\right) dt\right)$$
$$= -\lim_{\varepsilon \to 0} E\left(\int_{-\infty}^\infty p_\varepsilon(Y - t) \left(M_t \int_0^{S_T} \varphi(y) \frac{\partial R}{\partial t}(t, \tau_y) \, d\tau_y\right) dt\right).$$

Notice that

$$\left|\int_0^{S_T} \varphi(y) \frac{\partial R}{\partial s}(s, \tau_y) \, d\tau_y\right| \leq T \sup_{0 \leq s, u \leq T} \left|\frac{\partial R}{\partial s}(s, u)\right| \|\varphi\|_\infty.$$

Hence, applying the dominated convergence theorem, we get

$$E(M_Y) = 1 + \lambda E(\delta_t M|_{t=Y})$$
$$= 1 - \lambda \lim_{\varepsilon \to 0} E\left(\int_{-\infty}^{\infty} p_\varepsilon(Y-t)\left(M_t \int_0^{S_T} \varphi(y)\frac{\partial R}{\partial s}(t,\tau_y)\,d\tau_y\right)dt\right)$$
$$= 1 - \lambda E\left(M_Y \int_0^{S_T} \varphi(y)\frac{\partial R}{\partial s}(Y,\tau_y)\,d\tau_y\right). \qquad \square$$

The next step will be to replace the function $\varphi(x)$ by an approximation of the identity and let $T$ tend to infinity. Notice that (3.5) still holds for $\varphi(x) = \mathbf{1}_{[0,b]}(x)$ for any $b \geq 0$. In this way we can establish the following result.

THEOREM 3.4. *Assume conditions* (H1), (H2) *and* (H3). *For any* $a > 0$ *and* $\lambda \in \mathbb{R}$, *we have*

$$\int_0^a E(M_{\tau_y})\,dy$$
(3.9)
$$= a - \lambda E\left(\int_0^a \int_0^1 M_{z\tau_{y^+}+(1-z)\tau_y}\frac{\partial R}{\partial s}(z\tau_{y^+}+(1-z)\tau_y,\tau_y)\,dz\,d\tau_y\right).$$

Notice that we are not able to differentiate with respect to $a$, the integral in the rightmost expectation of (3.9), because the (random) measure $d\tau_y$, in general, is not absolutely continuous with respect to the Lebesgue measure.

PROOF OF THEOREM 3.4. Fix $a > 0$. We first replace the function $\varphi(x)$ by an approximation of the identity of the form $\varphi_\varepsilon(x) = \varepsilon^{-1}\mathbf{1}_{[0,1]}(x/\varepsilon)$ in formula (3.5). We will make use of the following notation:

$$Y_{\varepsilon,a} = \int_0^\infty \varphi_\varepsilon(x-a)(\tau_x \wedge T)\,dx.$$

At the same time we fix a nonnegative smooth function $\psi(x)$ with compact support such that $\int_\mathbb{R} \psi(a)\,da = c$ and we set

$$\int_\mathbb{R} E(M_{Y_{\varepsilon,a}})\psi(a)\,da$$
$$= c - \lambda \int_\mathbb{R} E\left(M_{Y_{\varepsilon,a}} \int_0^{S_T} \varphi_\varepsilon(y-a)\frac{\partial R}{\partial s}(Y_{\varepsilon,a},\tau_y)\,d\tau_y\right)\psi(a)\,da.$$

The increasing property of the function $x \to \tau_x$ implies that $\tau_{a^+} \wedge T \leq Y_{\varepsilon,a} \leq \tau_{a+\varepsilon} \wedge T$. Hence, $Y_\varepsilon$ converges to $\tau_{a^+} \wedge T$ as $\varepsilon$ tends to zero. Thus, almost surely, we have

$$\lim_{\varepsilon \to 0} M_{Y_{\varepsilon,a}} = \exp(\lambda X_{\tau_{a^+} \wedge T} - \tfrac{1}{2}\lambda^2 V_{\tau_{a^+} \wedge T}).$$

By the dominated convergence theorem,
$$\lim_{\varepsilon\to 0}\int_{\mathbb{R}} E(M_{Y_{\varepsilon,a}})\psi(a)\,da = \int_{\mathbb{R}} E(\exp(\lambda X_{\tau_{a+}\wedge T} - \tfrac{1}{2}\lambda^2 V_{\tau_{a+}\wedge T}))\psi(a)\,da.$$

Now, set $F(t) = M_t \frac{\partial R}{\partial s}(t,\tau_y)$. Then, assuming that $\varphi_\varepsilon(x) = \varepsilon^{-1}\mathbf{1}_{[0,1]}(x/\varepsilon)$, we have

$$\int_{y-\varepsilon}^{y} \varphi_\varepsilon(y-a) M_{Y_{\varepsilon,a}} \frac{\partial R}{\partial s}(Y_{\varepsilon,a},\tau_y)\psi(a)\,da$$
$$= \frac{1}{\varepsilon^2}\int_{y-\varepsilon}^{y}\mathbf{1}_{[0,1]}\left(\frac{y-a}{\varepsilon}\right)F\left(\int_{a}^{a+\varepsilon}\mathbf{1}_{[0,1]}\left(\frac{x-a}{\varepsilon}\right)(\tau_x\wedge T)\,dx\right)\psi(a)\,da$$
$$= \int_0^1 F\left(\int_0^1(\tau_{y+\varepsilon\xi-\varepsilon\eta}\wedge T)\,d\xi\right)\psi(y-\varepsilon\eta)\,d\eta$$
$$= \int_0^1 F\left(\int_0^\eta(\tau_{y+\varepsilon\xi-\varepsilon\eta}\wedge T)\,d\xi + \int_\eta^1(\tau_{y+\varepsilon\xi-\varepsilon\eta}\wedge T)\,d\xi\right)\psi(y-\varepsilon\eta)\,d\eta.$$

As $\varepsilon$ tends to zero, this expression clearly converges to

$$\psi(y)\int_0^1 F(\eta(\tau_y\wedge T) + (1-\eta)(\tau_{y+}\wedge T))\,d\eta.$$

So, we have proved that

(3.10)
$$\lim_{\varepsilon\to 0}\int_{\mathbb{R}} M_{Y_{\varepsilon,a}}\varphi_\varepsilon(y-a)\frac{\partial R}{\partial s}(Y_{\varepsilon,a},\tau_y)\psi(a)\,da$$
$$= \psi(y)\int_0^1 M_{z\tau_{y+}+(1-z)\tau_y}\frac{\partial R}{\partial s}(z\tau_{y+}+(1-z)\tau_y,\tau_y)\,dz.$$

In order to complete the proof of the theorem, we will apply the dominated convergence theorem. We have the following estimate for $y \leq S_T$:

$$\left|\int_{\mathbb{R}} M_{Y_{\varepsilon,a}}\varphi_\varepsilon(y-a)\frac{\partial R}{\partial s}(Y_{\varepsilon,a},\tau_y)\psi(a)\,da\right| \leq \|\psi\|_\infty \sup_{s,t\leq T}\left|\frac{\partial R}{\partial s}(s,t)\right|\sup_{t\leq T}|M_t|,$$

which allows us to commute the limit (3.10) with the integral with respect to the measure $P\times d\tau_y$ on the set $\{(\omega,y):y\leq S_T(\omega)\}$. In this way we get

$$\int_{\mathbb{R}} E(M_{\tau_y})\psi(y)\,dy$$
$$= \int_{\mathbb{R}}\psi(y)\,dy$$
$$\quad - \lambda E\left(\int_0^{S_T}\psi(y)\int_0^1 M_{z\tau_{y+}+(1-z)\tau_y}\frac{\partial R}{\partial s}(z\tau_{y+}+(1-z)\tau_y,\tau_y)\,dz\,d\tau_y\right).$$

Approximating $\mathbf{1}_{[0,a]}$ by a sequence of smooth functions $(\psi_n, n\geq 1)$ and letting $T$ tend to infinity completes the proof. $\square$

If we assume that the partial derivative $\frac{\partial R}{\partial t}(t,s)$ is nonnegative, then we can derive the following result.

PROPOSITION 3.5. *Assume that $X$ satisfies hypotheses* (H1), (H2) *and* (H3). *If* $\frac{\partial R}{\partial s}(s,t) \geq 0$, *then, for all $\alpha, a > 0$, we have*

(3.11) $$E(\exp(-\alpha V_{\tau_a})) \leq e^{-a\sqrt{2\alpha}}.$$

PROOF. Since $\frac{\partial R}{\partial t}(t,s) \geq 0$, we obtain
$$E(M_{\tau_a}) \leq 1,$$
that is,
$$E(\exp(\lambda a - \tfrac{1}{2}\lambda^2 V_{\tau_a})) \leq 1,$$
or
$$E(\exp(-\alpha V_{\tau_a})) \leq e^{-a\sqrt{2\alpha}}.$$
The result follows. □

The above proposition means that the Laplace transform of the random variable $V_{\tau_a}$ is dominated by the Laplace transform of $\tau_a$, where $\tau_a$ is the hitting time of the level $a$ for the ordinary Brownian motion. This domination implies some consequences on the moments of $V_{\tau_a}$. In fact, for any $r > 0$, we have, multiplying (3.11) by $\alpha^r$,

$$E(V_{\tau_a}^{-r}) = \frac{1}{\Gamma(r)} \int_0^\infty E(e^{-\alpha V_{\tau_a}}) \alpha^{r-1}\, d\alpha$$

(3.12) $$\leq \frac{1}{\Gamma(r)} \int_0^\infty e^{-a\sqrt{2\alpha}} \alpha^{r-1}\, d\alpha$$

$$= \frac{2^r \Gamma(r+1/2)}{\sqrt{\pi}} a^{-2r}.$$

On the other hand, for $0 < r < 1$,

$$E(V_{\tau_a}^r) = \frac{r}{\Gamma(1-r)} \int_0^\infty (1 - E(e^{-\alpha V_{\tau_a}})) \alpha^{-r-1}\, d\alpha$$

(3.13) $$\geq \frac{r}{\Gamma(1-r)} \int_0^\infty (1 - e^{-a\sqrt{2\alpha}}) \alpha^{-r-1}\, d\alpha.$$

In particular, for $r \in [1/2, 1)$, $E(V_{\tau_a}^r) = +\infty$.

REMARK 3.6. If $X$ is the standard Brownian motion, its covariance $s \wedge t$ does not satisfy condition (H1), but we still can apply our approach.

It is known from [3] that $d\tau_a$ has no absolutely continuous part and that $\{a, \tau_a = \tau_a^+\}$ is a Cantor set, hence, of zero Lebesgue measure. It follows from this observation and from (3.10) that

$$\int E(M_{\tau_y})\psi(y)\,dy = \int \psi(y)\,dy.$$

Choosing $\psi = \mathbf{1}_{[0,a]}$ yields to the expected result:

$$E\left(\int_0^a e^{\lambda y - (\lambda^2/2)V(\tau_y)}\,dy\right) = a.$$

If $X$ has independent increments and satisfies (H3), then

$$E(e^{-(\lambda^2/2)V(\tau_a)}) = e^{-\lambda a}.$$

This follows easily from the fact that $X$ can be written as a time-changed Brownian motion.

REMARK 3.7. Consider that $X$ is a fractional Brownian motion of Hurst index $H = 1$. Then $R(s,t) = st$, and consequently, $X_t = Yt$, where $Y$ is a one-dimensional standard Gaussian random variable. Then, $\tau_a = \tau_{a^+} = a/Y^+$. It is then easy to compute the Laplace transform of $\tau_a$ and we obtain

(3.14) $$E(\exp(-\alpha\tau_a^2)) = \tfrac{1}{2}e^{-a\sqrt{2\alpha}}.$$

We show now that our formula also yields to the right answer. We just note that $(y \mapsto \tau_y)$ is continuous. This entails that

$$\frac{\partial R}{\partial s}(z\tau_{y^+} + (1-z)\tau_y, \tau_y) = \frac{\partial R}{\partial s}(\tau_y, \tau_y) = \frac{1}{2}V'(\tau_y)$$

and

(3.15) $$\int_0^a E\left(\exp\left(\lambda y - \frac{\lambda^2}{2}V(\tau_y)\right)\right)dy$$
$$= a - \frac{\lambda}{2}E\left(\int_0^a \exp\left(\lambda y - \frac{\lambda^2}{2}V(\tau_y)\right)V'(\tau_y)\,d\tau_y\right).$$

Set

$$\Psi(a,\lambda) = E\left(\exp\left(\lambda a - \frac{\lambda^2}{2}V(\tau_a)\right)\right),$$

then

(3.16) $$\frac{\partial \Psi}{\partial a}(a,\lambda) = \lambda\Psi(a,\lambda) - \frac{\lambda^2}{2}E\left(M_{\tau_a}\frac{\partial V(\tau_a)}{\partial a}\right).$$

Substitute (3.15) into (3.16) to obtain

$$\frac{\partial \Psi}{\partial a} = 2\lambda\Psi - \lambda.$$

Then, there exists a function $C(\lambda)$ such that

$$\Psi(a,\lambda) = \frac{1}{2} + C(\lambda)e^{2\lambda a} \qquad \text{so that } E\left(\exp\left(-\frac{\lambda^2}{2}\tau_a^2\right)\right) = \frac{1}{2}e^{-\lambda a} + C(\lambda)e^{\lambda a}.$$

By dominated convergence, it is clear that, for any $\lambda$,

$$E\left(\exp\left(-\frac{\lambda^2}{2}\tau_a^2\right)\right) \stackrel{a\to\infty}{\longrightarrow} 0,$$

thus, $C(\lambda) = 0$ and

$$E\left(\exp\left(-\frac{\lambda^2}{2}\tau_a^2\right)\right) = \frac{1}{2}e^{-\lambda a}.$$

Changing $\lambda^2/2$ into $\alpha$ gives (3.14).

REMARK 3.8. Consider the case of a fractional Brownian motion with Hurst parameter $H > \frac{1}{2}$. Conditions (H1), (H2) and (H3) are satisfied and we obtain

$$\int_0^a E(M_{\tau_y})\,dy$$
$$= a - \lambda H E\bigg(\int_0^a \int_0^1 M_{z\tau_{y^+}+(1-z)\tau_y}([z\tau_{y^+} + (1-z)\tau_y]^{2H-1}$$
$$- |z(\tau_{y^+} - \tau_y)|^{2H-1})\,dz\,d\tau_y\bigg).$$

Moreover, $E(e^{-\alpha\tau_a^{2H}}) \leq e^{-a\sqrt{2\alpha}}$, and therefore, $E(\tau_a^p) < \infty$ if $p < H$. According to (3.13), $E(\tau_a^p)$ is infinite if $pH > 1/4$ and (3.12) entails that $\tau_a$ has finite negative moments of all orders.

**Acknowledgments.** This work was carried out during a stay of Laurent Decreusefond at Kansas University, Lawrence, KS. He would like to thank KU for warm hospitality and generous support.

Department of Mathematics of Information
Communications and Calculus
Telecom Paris
46, rue Barrault
75634 Paris
France
E-mail: Laurent.Decreusefond@enst.fr

Department of Mathematics
University of Kansas
405 Snow Hall
1460 Jayhawk Blvd
Lawrence, Kansas 66045-7523
USA
E-mail: nualart@math.ku.edu